\documentclass{article}

\newtheorem{theorem}{Theorem}[section]
\newtheorem{lemma}[theorem]{Lemma}
\newtheorem{corollary}[theorem]{Corollary}
\newtheorem{definition}[theorem]{Definition}

\newcommand\prob{{\mbox Prob}}
\newcommand\qed{\begin{flushright} {\bf q.e.d.} \end{flushright} }
\newcommand\prf{\noindent {\bf Proof :}  }

\newcommand\np{\mbox{NP}}

\newcommand\nn{{\{0,1\}^n}}

\newcommand\s{ \{0, 1\}^s }

\newcommand\NN{{\{0,1\}^{N}}}

\newcommand\at{ \{0, 1\}^t }
\newcommand\ar{ \{0, 1\}^r }
\newcommand\as{ \{0, 1\}^s }
\newcommand\au{ (u, q^u, r^v) }
\newcommand\av{ (v, q^u, r^v) }
\newcommand\bits{\{0,1\}}

\newcommand\ff{{{\bf F}_2}}
\newcommand\rlin{{\mbox{R(LIN}/{\bf F}_2)}}
\newcommand\newax{Ax(\pi, {\bf r})}
\newcommand\lw{\mbox{lw}}
\newcommand\rr{{\bf r}}
\newcommand\bs{{\bf s}}
\newcommand\pr{{\pi}({\rr})}

\newcommand\cc{\mbox{CC}}
\newcommand\pcc{ C^{pub}_{\epsilon} }
\newcommand\mc{\mbox{MCC}_U}
\newcommand\ccr{\mbox{CC}^{{\bf R}} }

\newcommand\pro{(G, {\sf lab}, F, S)}
\newcommand\pp{{\mathbf P}}
\newcommand\rpp{{\mathbf P}_{\bf r}}
\newcommand\rrpp{({\mathbf P}_{\bf r})_{\bf r}}
\newcommand\rpro{(G_{\bf r}, {\sf lab}_{\bf r}, F_{\bf r}, S_{\bf r})}

\newcommand\trrpp{(\tilde{\mathbf P}_{\bf r})_{\bf r}}
\newcommand\rcc{{C}_{\bf r}}
\newcommand\rrcc{({C}_{\bf r})_{\bf r}}

\newcommand\clique{{Clique_{n_0,\omega}}}
\newcommand\color{{Color_{n_0,\xi}}}

\newcommand\umin{{U^{min}}}
\newcommand\vmax{{V^{max}}}

\begin{document}

\title{Randomized feasible interpolation and
monotone circuits with a local oracle}

\author{Jan Kraj\'{\i}\v{c}ek}

\date{Faculty of Mathematics and Physics\\
Charles University in Prague}

\maketitle

\begin{abstract}

The feasible interpolation theorem for semantic derivations from K. (1997) \cite{Kra-interpol} 
allows to derive from some short semantic
derivations (e.g. in resolution)
of the disjointness of two $\np$ sets $U$ and $V$
a small communication protocol (a general dag-like
protocol in the sense of K. (1997) \cite{Kra-interpol})
computing the Karchmer-Wigderson multi-function $KW[U,V]$ associated with the sets, and such a protocol 
further yields a small circuit separating $U$ from $V$. When $U$ is closed upwards the protocol
computes the monotone Karchmer-Wigderson multi-function $KW^m[U,V]$ and the resulting circuit
is monotone.
K. (1998) \cite{Kra-game}
extended the feasible interpolation theorem
to a larger class of semantic derivations
using the notion of a real communication complexity (e.g. to the cutting planes
proof system CP).

In this paper we generalize the method to a still larger class of
semantic derivations by allowing randomized protocols.  
We also introduce an extension of the monotone circuit model,
monotone circuits with a local oracle (CLOs), that does correspond to communication protocols 
for $KW^m[U,V]$ making errors. The new randomized feasible interpolation thus shows that a short
semantic derivation (from a certain class of derivations larger than in the original method) 
of the disjointness of $U, V$, $U$ closed upwards, yields a small randomized protocol for
$KW^m[U,V]$ and hence a small monotone CLO
separating the two sets. 

This research is motivated by the open problem to establish a lower bound for
proof system $\rlin$ operating with clauses formed by linear Boolean functions over $\ff$.
The new randomized feasible interpolation applies to this proof system and also 
to (the semantic versions of) cutting planes CP, to small width resolution
over CP  of K. (1998) \cite{Kra-modules} (system R(CP)) and to random resolution RR
of Buss, Kolodziejczyk and Thapen \cite{BKT}. The method does not yield yet
lengths-of-proofs lower bounds; for this it is necessary to establish lower bounds for
randomized protocols or for monotone CLOs.

\end{abstract}

Consider a propositional proof system 
$\rlin$
that operates with
clauses of linear equations over $\ff$ and combines
the rules of both resolution and linear equational calculus.
A line $C$ in a proof has the form
$$
\{f_1, \dots, f_k\}
$$
with $f_i \in \ff[x_1, \dots, x_n]$ linear polynomials
and the intended meaning is that an assignment $x := a \in \nn$
to variables makes
$C$ true if and only if
one of $f_i = 1$ becomes true, i.e. the truth value of $C$
is computed by Boolean formula
$$
\bigvee_{i \le k} f_i
$$
in the language with $\bigvee, \oplus, 0, 1$. We often leave 
the outside
brackets $\{, \}$ out when writing clauses. For $L \subseteq C$ define
$\sum L := \sum_{f \in L} f$.

The rules of $\rlin$ are the following four:
$$
\frac{}{h, h+1}\ \ \ \ \ \ \frac{C}{C,f}\ \ \ \ \ \ \frac{C,0}{C}
\ \ \ \ \ \ 
\frac{C,g\ \ \ \ C,h}{C,g+h+1}\ .
$$
We shall call the rules {\em $\ff$-axiom}, {\em weakening}, {\em contraction} and 
{\em the binary rule}, respectively.
This proof system (albeit defined slightly differently but polynomially equivalently, denoted Res-Lin 
there) has been
considered already by Itsykson and Sokolov \cite{ItsSok14} 
who proved an exponential lower bound
for tree-like proofs. They also showed that 
the semantic version of the system
(in the sense of semantic derivations of \cite{Kra-interpol}) is p-equivalent
to the syntactic version, whether tree-like or dag-like.
This paper is motivated by the problem to establish a lower bound for unrestricted (i.e. dag-like) $\rlin$ proofs.

Proof systems combining resolution or, more generally, logical reasoning 
with algebraic reasoning were considered earlier
by several authors: \cite{Kra-modules}
defined proof systems $R(CP)$ and $LK(CP)$ extending cutting plane by
a logic reasoning and proved an exponential lower bound for a subsystem
of $R(CP)$, 
Hirsch and Kojevnikov \cite{HirKoj,Koj} considered resolution over 
a system for linear programing and Kojevnikov \cite{Koj} improved upon a bound in
\cite{Kra-modules}.
Raz and Tzameret \cite{RazTza08} studied resolution over linear
equations with integral coefficients
and proved a lower bound for a class of its proofs, 
and Alekhnovich et.al. \cite{ABRW} defined polynomial calculus with resolution
PCR which extends PC in a way that incorporates resolution (lines of
proofs are polynomials, however).

There is also a link to the well-known open problem to establish
lower bounds for constant depth Frege systems in DeMorgan language
augmented by a connective counting modulo a prime, the so called
$AC^0[p]$-Frege systems. The strongest subsystem 
of such a system for which a lower bound is
known is a low degree polynomial calculus
operating with polynomials formed from $AC^0$-formulas, \cite{Kra-mfcs}.
The lower bound problem for $\rlin$ seems interesting also
because the top proof system is logical. 
Note that Buss, Kolodziejczyk and Zdanowski \cite{BKZ} proved that, in fact,
the $AC^0[p]$-Frege system collapses (with a quasi-polynomial blow-up
in proof size) to a proof system operating with clauses of conjunctions of low
degree polynomials.

Our approach is to use feasible interpolation for semantic derivations
from \cite{Kra-interpol} but we
need to generalize it first to allow small errors. 
The generalization we develop here allows randomized communication protocols with errors
(protocols in the sense of \cite{Kra-interpol}) for computing the Karchmer-Wigderson multi-function. 
Protocols making no errors correspond to separating circuits but protocols with errors do
not yield separating circuits making some error. Instead we introduce an extension of the circuit model,
circuits with a local oracle (CLO), that does correspond to protocols with errors.

Tree-like protocols with errors for $KW^m[U,V]$
yield monotone separating formulas with a local oracle and subsume
the ordinary Karchmer-Wigderson (1988) \cite{KW} protocols pictured as binary trees. 
A lower bound in this case is known (cf. \cite{IPU,Kra-game}
for examples based on the bipartite perfect matching problem and Hall's theorem).
Further, monotone CLOs efficiently simulate monotone real circuits (Section \ref{20.3.16c})
and any two disjoint sets can be 
separated by a small non-monotone CLO (Lemma \ref{24.11.16a}
and the remark at the end of Section \ref{24.9.16a}).
To establish a lower bound for monotone CLOs separating two
$\np$ sets, one closed upwards, is an open problem.

To be able to apply randomized feasible interpolation to $\rlin$ we use 
the approximation method of Razborov \cite{Raz87} and Smolensky \cite{Smo}
in order to reduce the linear width 
(defined in Section \ref{4.9.15b})
in a general not too long proof at the expense 
of introducing an error (cf. Section \ref{4.9.15c}). 
The new method may have further applications and, in particular, it applies to 
the semantic versions of cutting planes CP, to small width resolution over cutting planes 
R(CP), and to random resolution RR.
The method on its own does not yield yet
lengths-of-proofs lower bounds; for this it is necessary to establish lower bounds for
randomized protocols or for monotone CLOs. Some partial results about monotone CLOs are 
obtained in \cite{KraOli}.

\medskip

The paper is organized as follows.
Section \ref{20.3.16a} recalls some notions and
results from \cite{Kra-interpol}.
In Section \ref{20.3.16b} 
we define the concept of randomized protocols and use it to formulate randomized feasible interpolation.
In Section \ref{24.9.16a} we introduce circuits with a local oracle (CLO)  and prove that they correspond to
protocols with errors and that, in particular, randomized protocols yield CLOs.
In Section \ref{4.9.15b} we introduce the linear width of $\rlin$ proofs and discuss 
the case when it is small. 
Randomized feasible interpolation is proved for $\rlin$ in Section \ref{4.9.15c}
and for CP and small width R(CP) in Section \ref{20.3.16c}.
The lower bound problem for monotone CLOs (and hence for randomized
protocols computing the monotone Karchmer-Wigderson multi-function for some pair of sets) 
is discussed in Section \ref{20.3.16d}.
The paper is concluded by a few remarks in Section \ref{4.9.15e}.
A proof complexity background can be found in \cite{kniha,Pud-survey}.

\section{Feasible interpolation preliminaries} \label{20.3.16a}

The general feasible interpolation theorem from \cite{Kra-interpol} 
for semantic derivations uses communication complexity. 
One considers two disjoint NP sets $U, V \subseteq \nn$ and the Karchmer-Wigderson 
multi-function whose valid values on a pair $(u,v) \in U \times V$ is any
coordinate in which $u, v$ differ.
The aim is to extract from a short proof of the disjointness of $U, V$ 
some upper bound on the 
computational complexity of this multi-function in some computational 
model. Proving then a computational complexity lower bound for the model
allows to infer
a length-of-proofs lower bound. The original set-up (and the one most 
frequently used) derives from 
the proof data the existence of a small circuit
separating $U$ and $V$. In the monotone case one can use then known strong 
lower bounds for monotone circuits, for example Alon and Boppana \cite{AloBop}.

When the construction of \cite{Kra-interpol} is applied to tree-like proofs 
it leads to familiar 
protocols for communication that are pictured as binary trees, cf.\cite{KW}. 
However, for applications to general, dag-like, proofs one needs 
a more general notion of a protocol defined in \cite[Def.2.2]{Kra-interpol}. 
The key fact, allowing to prove some lower bounds, is that similarly as small tree-like communication protocols 
correspond to small formulas separating $U$ and $V$ (by Karchmer and Wigderson \cite{KW}), 
the more general protocols used in \cite{Kra-interpol}
correspond to small separating circuits.

Let us now recall formally relevant definitions and facts from \cite{Kra-interpol}.
A multi-function on $U \times V$ with values in some set $I \neq \emptyset$
is a ternary relation 
$R \subseteq U \times V \times I$ such that for all $(u,v) \in U \times V$
there is $i \in I$ such that $R(u,v,i)$.
Some value for $(u,v)$ from its domain can be computed
by two players, one receiving $u$ and the other one $v$,
exchanging bits of information until they agree on a valid value $i$. 
The communication complexity of $R$, $CC(R)$, is the minimal number of bits
they need to exchange (in an optimal protocol) in the worst case.

The Karchmer-Wigderson multi-function $KW[U,V]$ of a particular interest is defined 
for two disjoint sets $U,V \subseteq \nn$: a valid value 
of $KW(u,v)$ on pair $(u,v) \in U \times V$
is any $i \in [n]$ such that $u_i \neq v_i$.
The monotone version of this function $KW^m[U,V]$ is defined when $U$ 
is closed upwards (or $V$ downwards) 
and a valid value on $(u,v)$ is any $i \in [n]$ such that
$u_i =1 \wedge v_i =0$.

Given two disjoint $U, V \subseteq \{0, 1\}^n$ 
and $R \subseteq U \times V \times I$ a multi-function, \cite[Def.2.2]{Kra-interpol}
defines a protocol for $R$ to be a 4-tuple $\pp = \pro$
satisfying the following conditions:
\begin{enumerate}

\item[(P1)] $G$ is a directed acyclic graph that has one source (the in-degree $0$ node
called the {\em root}) denoted $\emptyset$.

\item[(P2)] The nodes with the out-degree $0$ are {\em leaves}
and they are labelled by the mapping ${\sf lab}$ by
elements of $I$.

\item[(P3)] $S(u,v,x)$ is a function (the {\em strategy})
that assigns to a node $x \in G$ 
and a pair $u \in U$ and $v \in V$
node $S(u,v,x)$ accessible by an edge from $x$.

\item[(P4)] For every 
$u \in U $ and $v \in V$, $F(u,v) \subseteq
G$ is a set (called the {\em consistency condition}) satisfying:
\begin{enumerate}

\item $\emptyset \in F(u,v)$,

\item $x \in F(u,v) \rightarrow S(u,v,x) \in F(u,v)$,

\item if $x \in F(u,v)$ is a leaf and ${\sf lab}(x)=i$, then
$R(u,v,i)$ holds.
\end{enumerate}
\end{enumerate}

\noindent
We say that $\pp$ is {\em tree-like} iff $G$ is a
tree. 

\medskip

The complexity of $\pp$ is measured by its {\em size}, which is the cardinality of $G$, 
and by the following notion: The {\em communication complexity} of $\pp$, denoted
$CC(\pp)$, is the minimal $t$ such that  for every $x \in G$
the communication complexity for 
the players (one knowing $u$ and $x$, the other one $v$ and $x$)
to decide $x \in_? F(u,v)$ or to compute $S(u,v,x)$
is at most $t$.

\bigskip

The interpolation theorem in \cite{Kra-interpol} was formulated using the
notion of a {\em semantic derivation} (\cite[Def.\ 4.1]{Kra-interpol}): 
A sequence of sets $D_1, \dots, D_k \subseteq \{0,1\}^N$ 
is a semantic derivation of $D_k$ from $A_1, \dots, A_m \subseteq \NN$
if each $D_i$ is either one of $A_j$'s or contains
$D_{j_1} \cap D_{j_2}$, for some $j_1, j_2 < i$. A semantic derivation is a {\em refutation} 
of $A_1, \dots, A_m$ iff $D_k = \emptyset$.

We shall introduce now a general set-up for our investigation of interpolation and
 we shall refer to it the whole paper. We assume the following
conditions for parameters and sets, and introduce the following notation:

\begin{equation}
N = n+s+r\ ,\ N, n \geq 1\ .
\end{equation}

\begin{equation}
A_1, \dots, A_m \subseteq \{0, 1\}^{n+s}\ \mbox{ and }\ 
B_1, \dots, B_{\ell} \subseteq \{0, 1\}^{n+r}\ .
\end{equation} 
From the total $N$ variables, $n$ represent an input $a$ from $\nn$, $s$ variables represent a potential witness
$b$ for the membership of $a$ in $U$ and $r$ variables represent a potential witness
$c$ for the membership of $a$ in $V$ ($U$ and $V$ are defined below).
For $A \subseteq \{0, 1\}^{n+s}$ define
\begin{equation}
\tilde A := \bigcup_{(a,b) \in A} \{(a,b,c)\ |\ c \in \ar\}
\end{equation}
and for $B \subseteq \{0, 1\}^{n+r}$ define:
\begin{equation}
\tilde B := \bigcup_{(a,c) \in B} \{(a,b,c)\ |\ b \in \as\}\ .
\end{equation}
where $a, b, c$ range over $\nn$, $\as$ and $\ar$, respectively.
Define:
\begin{equation}
U = \{ u \in \nn \ |\ \exists b \in \as; 
(u, b) \in \bigcap_{j \le m} A_j\}
\end{equation}
and
\begin{equation}
V = \{ v \in \nn \ |\ \exists c \in \ar; 
(v, c) \in \bigcap_{j \le \ell} B_j\}\ .
\end{equation}
We shall also refer to the following monotonicity condition.
For all $u, u' \in \nn$ and $b \in \s$:
\begin{equation}
(u,b) \in \bigcap_{j \le m} A_j \wedge u' \geq u 
\longrightarrow (u',b) \in \bigcap_{j \le m} A_j\ .
\end{equation}

\medskip

The complexity of sets in a semantic derivation is measured by
the following notion of (monotone) communication complexity 
of subsets of $\NN$ defined in
\cite{Kra-interpol}.
For $D \subseteq \NN$,
$u, v \in \nn$, $q^u \in \as$ and $r^v \in \ar$
consider four tasks:
\begin{enumerate}

\item Decide whether $\au \in D$.

\item Decide whether $\av \in D$.

\item  If $\au \in D \not\equiv \av \in D$ find $i \le n$
such that $u_i \neq v_i$.

\item If $\au \in D$ and $\av \notin D$ either find $i \le n$
such that 
$$
u_i =1 \wedge v_i =0
$$
or decide that there is some $u'$ satisfying 
$$
u'\geq u \wedge (u', q^u, r^v) \notin D\ .
$$
\end{enumerate}
The {\em communication complexity} $CC(D)$ of $D$ 
is the minimal $t$ such that the tasks 1.-3.
can be solved by the players, one knowing
$u, q^u$ and the other one knowing $v, r^v$, exchanging 
at most $t$ bits.
The {\em monotone communication complexity w.r.t. $U$} of $D$, denoted $\mc(D)$,
is the minimal $t \geq CC(D)$ such that also the task 4.
can be solved by the players exchanging 
at most $t$ bits.

Now we are ready to recall a fact about the existence of
protocols from the proof of 
\cite[Thm.5.1]{Kra-interpol}.

\begin{theorem}[{\cite{Kra-interpol}}] \label{3.2}
{\ }

Assume the set-up conditions (1)-(6) and assume that 
$\pi = D_1, \dots, D_k$ is a semantic refutation of 
the sets $\tilde A_1, \dots, \tilde A_m, 
\tilde B_1, \dots, \tilde B_{\ell}$. Let $t \geq 1$ be such that $t \geq CC(D_i)$
for all $i \le k$.

Then there is a protocol for $KW[U,V]$ of size $k + 2n$ and of communication complexity 
$O(t)$.
The protocol has $k$ inner vertices, the sets in $\pi$, and additional $2n$  vertices, the leaves,\
labelled by all possible formulas $u_i = 1 \wedge v_i = 0$ and $u_i = 0 \wedge v_i = 1$.

If condition (7) is also satisfied and $\mc(D_i) \le t$ for all $i \le k$ then there
is a protocol for $KW^m[U,V]$ of size $k+n$ and of communication complexity
$O(t)$.

Further, the consistency condition $F$ is defined in both the monotone and the non-monotone cases
identically as:
$$
D \in F(u,v)\ \mbox{ iff }\ \av \notin D
$$
for $D$ in $\pi$, and 
$$
x \in F(u,v)\ \mbox{ iff }\ {\sf lab}(x) \mbox{ is valid for } u, v
$$
for $x$ a leaf.

Moreover, if $\pi$ is tree-like, so is $G$.
\end{theorem}

\section{Randomized feasible interpolation for semantic derivations} \label{20.3.16b}

First we generalize protocols to allow a randomization and some error.

\begin{definition} \label{24.9.16c}
A {\em randomized protocol} for multi-function $R \subseteq U \times V \times I$
with error $\epsilon > 0$ is a random variable $\rrpp$ where each $\rpp$ is
a 4-tuple satisfying conditions (P1), (P2), (P3) and (P4a) defining
protocols and instead of conditions (P4b) and (P4c) it satisfies:

\begin{enumerate}
\item[(P4b')] For every $(u,v) \in U \times V$,
$$
\prob_{\bf r}[\exists x,\ x \in F_{\bf r}(u,v) \wedge
S_{\bf r}(u,v,x) \notin F_{\bf r}(u,v)] \ \le \ \epsilon\ .
$$

\item[(P4c')] For every $(u,v) \in U \times V$,
$$
\prob_{\bf r}[\exists \mbox{leaf } x,\ x \in F_{\bf r}(u,v) \wedge
{\sf lab}_{\bf r}(x) = i \wedge \neg R(u,v,i)] \ \le \ \epsilon\ .
$$
\end{enumerate}
The size of $\rrpp$ is $\max_{\bf r} size(\rpp)$ and
the communication complexity of $\rrpp$ is 
$\max_{\bf r} CC(\rpp)$. We say that $\rrpp$ is tree-like if each $\rpp$ is.

\end{definition}

We note a simple observation.

\begin{lemma} \label{24.9.16b}
For any randomized protocol $\rrpp$ 
for multi-function $R \subseteq U \times V \times I$
of size $S$, communication complexity $t$
and error $\epsilon$ there exists
a randomized protocol $\trrpp$ 
for multi-function $R$
of size at most $2S$ (with at most $S$ leaves), communication complexity
at most $3t$ and error $\epsilon$
such that $(P4b)$ never fails, i.e. the probability in $(P4b')$ 
is $0$.

\end{lemma}

\prf

Introduce for each inner node $x \in G_{\bf r}$ a new leaf node $\tilde x$,
label it arbitrarily (e.g. $u_1 = 1 \wedge v-1 =0 $),
and define a new strategy $\tilde S_{\rr}$ that first checks if 
$$
x \in F_{\bf r}(u,v) \rightarrow S_{\bf r}(u,v,x)\in F_{\bf r}(u,v)
$$
is true and if so it uses $S_{\bf r}$, otherwise it sends $x$ into $\tilde x$
and the failure of the condition is the definition of $\tilde x\in 
{\tilde F}_{\bf r}(u,v)$.

\qed

In connections with interpolation we are interested in the situation when the multi-function
is the Karchmer-Wigderson one. It makes sense to consider only the monotone case $KW^m[U,V]$
as the next lemma recalls.

\begin{lemma} [Raz and Wigderson \cite{RazWig89}] \label{24.11.16a}
Let $U, V$ be any two disjoint subsets of $\nn$. Then for any $\epsilon > 0$
there is a tree-like randomized protocol
$\rrpp$ computing $KW[U,V]$ of size $S = (n+ \epsilon^{-1})^{O(1)}$, 
communication complexity $t = O(\log n + \log(\epsilon^{-1}))$ and error $\epsilon$.

In particular, for $\epsilon = n^{-\Omega(1)}$ the size is $S = n^{O(1)} $
and the communication complexity is $t = O(\log n)$.

\end{lemma}

\prf

A randomized protocol computing $KW[U,V]$ 
is determined by $\log (\epsilon^{-1})$ subsets 
$I \subseteq [n]$. The players exchange the parity of the bits in their respective 
strings belonging to the first
such $I$, then to the second, etc. until they find $I$
for which the parity differs. Then 
they find a valid value for $KW[U,V]$ by binary search. If they do not find such $I$, they declare an error.
This gives a randomized protocol of size polynomial in 
$n, \epsilon^{-1}$, with communication complexity $2 (\log n + \log (\epsilon^{-1}))$, and error
$\epsilon$.

\qed

Now we introduce a notion that we will use in the context of semantic
derivations. Let $X \in \NN$ and let ${\cal Y} = (Y_{\bf r})_{\bf r}$
be a random distribution on subsets of $\NN$, and let $\delta > 0$.
We say that ${\cal Y}$ is a $\delta$-{\em approximation}
of $X$ iff for all $ w \in \NN$:
$$
\prob_{\bf r}[w \in X \triangle Y_{\bf r}]\ \le\ \delta\ 
$$
where $X \triangle Y$ is the symmetric difference.

Working in the set-up (1)-(6) the sets $X$ and $Y_{\bf r}$ are subsets of $\NN$ and the definitions of 
$CC$ and $MCC_U$ apply to them. With this in mind
we further define that the {\em (monotone) communication complexity} of ${\cal Y}$
is at most $t$ if this is true for all $Y_{\bf r}$, and that
the {\em $\delta$-approximate (monotone) communication complexity} of $X$ is at most
$t$ if there is a $\delta$-approximation ${\cal Y}$
of $X$ with this property.

\begin{theorem} \label{29.3.16b}
{\ }

Assume the set-up conditions (1)-(7). 
Let $\pi\ :\ D_1, \dots, D_k = \emptyset$ be a semantic refutation of
sets $\tilde A_1, \dots, \tilde A_m, \tilde B_1, \dots, \tilde B_{\ell}$ 
such that the $\delta$-approximate monotone
communication complexity of every $D_i$ is at most $t$.

Then  there is a randomized protocol $\rrpp$ for $KW^m[U,V]$ of size at most $k+n$,
communication complexity $O(t)$ and of error
at most $3 \delta k$.

Moreover, if the refutation $\pi$ is tree-like
then also $\rrpp$ is tree-like.

\end{theorem}

\prf

Take the protocol $\pp = \pro$ provided by Theorem \ref{3.2}. Its strategy $S$
and the consistency condition $F$ are defined in terms of sets $D_i$. In particular,
for any $(u,v)\in U \times V$ and $x \in G$ an inner node, both the value of $S(u,v,x)$ and
the truth value of $x \in F(u,v)$ are defined from at most $3$ truth values of
statements of the form $(u,q^u,r^v) \in D_i$ or $(v,q^u, r^v) \in D_i$ for some 
specific indices $i \le k$ determined by $x$, where $q^u$ and
$r^v$ depend just on $u$ and $v$, respectively.

Not knowing anything about the monotone communication complexity of the sets 
$D_i$ we cannot estimate the communication complexity of $\pp$.
At this point we use the $\delta$-approximations of the sets $D_i$. If 
$(E^i_{\bs})_{\bs}$ are $\delta$-approximations of $D_i$, $i \le k$, 
let the space of samples $\rr$ for $\rpp$ be the product of the sample spaces
of these $k$ $\delta$-approximations and define $S_{\rr}$ and $F_{\rr}$
as $S$ and $F$ before but using the particular sets $E^i_{\bs}$ 
(with $\bs$ determined by $\rr$) in place of the 
sets $D_i$. In particular, $D_i \in F_{\bf r}(u,v)$ iff $\av \notin E^i_{\bf s}$.
Further, put $G_{\rr} := G$ and ${\sf lab}_{\rr}:={\sf lab}$.

For any given $(u,v) \in U \times V$ and $x \in G$ the (truth) value of $S_{\rr}$
and $F_{\rr}$ differs from $S$ and $F$ respectively with probability at most $3 \delta$.
Hence for $(u,v)$ the error in conditions (P4b') and (P4c ') is at most 
$\epsilon := 3 \delta k$.

\qed

We describe yet another type of semantic refutations that also yields
randomized protocols.

\begin{theorem} \label{29.3.16a}
Assume the set-up conditions (1)-(7).
Let $e \geq 1$, $\epsilon > 0$ and let
$(C_{\bf r})_{\bf r}$ be a random distribution on 
$({\cal P}(\NN))^e$, $C_{\bf r} = (C^1_{\bf r}, \dots, C^e_{\bf r})$,
such that each $\bigwedge_{i \le e} (C^i_{\bf r})_{\bf r}$ is 
an $\epsilon$-approximation of $\NN$.

Assume that for all samples $\bf r$ there is a semantic refutation $\pi_{\bf r}$
of 
$$
\tilde A_1, \dots, \tilde A_m, \tilde B_1,\dots,\tilde B_{\ell}, C^1_{\bf r}, \dots, C^e_{\bf r}
$$
with $k$ lines, and such that the monotone communication complexity of all
sets in $\pi_{\bf r}$ is at most $t$.

Then there is a randomized protocol for $KW^m[U,V]$) of size at most 
$k+n+e \le 2k +n$,
communication complexity $O(t)$ and of error
at most $\epsilon$.

Moreover, if the refutation $\pi$ is tree-like
then also $\rrpp$ is tree-like.

\end{theorem}

\prf

The construction of $\pp = \pro$ in Theorem \ref{3.2} yields $G$ whose
inner nodes correspond to lines of the refutation and leaves are extra $n$ nodes. 
The construction uses the fact that
for $(u,v) \in U \times V$ the strings $q^u$ and $r^v$
are chosen so that $(u,q^u, r^v)\in\tilde A_i$ and $\av \in \tilde B_j$
for all $i \le m$ and $j\le \ell$.
In particular, each initial set $\tilde A_1, \dots, \tilde A_m, \tilde B_1, \dots, \tilde B_\ell$
contains either $\au$ or $\av$.

In the presence of the new initial clauses
$C^i_{\rr}$ this is no longer true and it may happen
that both $(u,q^u, r^v)$ and $(v, q^u, r^v)$ are outside of some
$C^i_{\rr}$.

We define $\rpp$ as follows. Each $G_{\rr}$ has $e$ extra leaves $y_i$
labelled arbitrarily (say $u_1 = 1 \wedge v_1=0$ for the definiteness)
and the strategy $S_{\rr}(u,v,x)$ sends node $x$ corresponding to 
$C^i_{\rr}$ to $y_i$ if
\begin{equation} \label{28.3.16a}
\av \notin C^i_{\rr}
\end{equation}
and the same condition defines when $y_i \in F_{\rr}(u,v)$.

As $\bigwedge_{i \le e} (C^i_{\rr})_{\rr}$ is an $\epsilon$-approximation of $\NN$,
(\ref{28.3.16a}) happens with probability at most $\epsilon$
in total.

\qed

\section{Monotone circuits with a local oracle}
\label{24.9.16a}

Our aim in this section is to define a generalization of the circuit model
that corresponds to protocols with errors computing $KW^m[U,V]$. 
We restrict ourselves to the monotone case due to Lemma \ref{24.11.16a}
(see also the remark at the end of this section).

A monotone {\em circuit with a local oracle} (monotone CLO, briefly) separating $U$ from $V$
is determined by the following data:

\begin{enumerate}

\item a monotone Boolean circuit $D(x_1,\dots, x_n, y_1, \dots, y_e)$ with inputs
$\overline x$ and $\overline y$,

\item a set $\cal R$ of combinatorial rectangles $U_j \times V_j \subseteq U \times V$, for $j \le e$,
called {\em oracle rectangles} of the CLO,

\end{enumerate}
and satisfying the following condition:

\begin{enumerate}

\item[3.] for all monotone Boolean functions $f_j : \nn \rightarrow \bits$, $j \le e$,
such that
$$
f_j(U_j) \subseteq \{1\}\ \ \ \ \mbox{ and }\ \ \ \ f_j(V_j) \subseteq \{0\}
$$
the function
$$
C(\overline x)\ :=\ D(\overline x, f_1(\overline x), \dots, f_e(\overline x))
$$
separates $U$ from $V$:
$$
C(U) = \{1\}\ \ \ \ \mbox{ and }\ \ \ \ C(V) = \{0\}\ .
$$
\end{enumerate}
The {\em size} of the CLO is the size of $D$ and its {\em locality} is
$$
\frac{|\bigcup_{j \le e} U_j \times V_j|}{|U \times V|}
$$
(we assume both $U, V$ are non-empty).
Note that $C$ defines a monotone Boolean function for any choice of monotone functions $f_j$.

\bigskip

The proof of the following lemma expands a bit upon a proof by Razborov \cite{Raz95}.

\begin{lemma} \label{25.9.16b}
Assume that $\rrpp$ is a randomized protocol for $KW^m[U,V]$ of size $s$, communication
complexity $t$ and error $\epsilon$.

Then there is a monotone circuit with a local oracle separating $U$ from $V$ 
of size $s 2^{O(t)}$ and locality $\epsilon$.
\end{lemma}

\prf

Assume $\rrpp$ is a randomized protocol 
satisfying the hypothesis of the lemma, with $\rpp = \rpro$. By 
Lemma \ref{24.9.16b}
we may assume that each $G_{\bf r}$ makes errors only in leaves, i.e. violates possibly only the
condition (P4c) of Section \ref{20.3.16a} in the sense of (P4c') of Definition \ref{24.9.16c}.
This may increase the size and the communication complexity proportionally but that does not
change the form $s 2^{O(t)}$ of the upper bound.

By averaging there must be some sample $\bf r$ such that $\rpp$ makes an error for at most 
$\epsilon$-part of all pairs $U \times V$. Fix one such protocol $\pro := \rpp$ for the rest of the proof.
We may also assume that the communication of the players deciding that a leaf $a$ is in  $F(u,v)$ 
ends with each player sending the value of the $i$-th bit of $u$ or $v$, respectively, where
$i = {\sf lab}(a)$. That is, they both know at the end whether an error occurred for $(u,v)$
and the set of these erroneous pairs for which $a \in F(u,v)$ is a disjoint union of
combinatorial rectangles.

For a vertex $a$ of $G$ and 
a string $w \in \at$ denote:

\begin{itemize}

\item
$R_{a,w}$ the rectangle 
$U_{a,w} \times V_{a,w}$, some $U_{a,w} \subseteq U$ and
$V_{a,w} \subseteq V$, 
of pairs $(u,v) \in U \times V$
such that the communication of the players
deciding $a \in_? F(u,v)$ evolves according to $w$
and ends with the affirmation of the membership,

\item
 $k_a$: the number of nodes in $G$
that can be reached from node $a$ by a directed path (so
$k_a = 1$ for $a$ a leaf, while $k_{\emptyset} \le s$
for the root $\emptyset$).

\end{itemize}
Assume 
\begin{equation} \label{25.9.16a}
R_1 := U_1 \times V_1, \dots, R_e := U_e \times V_e\ ,\ \mbox{ for } j \le e
\end{equation}
enumerate all rectangles $R_{a,w}$ where $a$ is a leaf and $(u,v) \in R_{a,w}$
iff $a \in F(u,v)$ and the players decided this with communication $w$ 
but ${\sf lab}(a)$ is incorrect for $(u,v)$, i.e. an error occurs
for $(u,v)$ at $a$.

\medskip
\noindent
{\bf Claim 1:} {\em
For all $a \in G$ and $w \in \at$ there is a size 
$\le k_a 2^{O(t)}$ monotone circuit with a local oracle
separating $U_{a,w}$ 
from $V_{a,w}$ such that its oracle rectangles are included among (\ref{25.9.16a}). 
The constant implicit in the exponent $O(t)$ is independent of $a$.
}

\smallskip
\noindent

For $a,w$ we shall denote by $D_{a,w}, {\cal R}_{a,w}$ a monotone CLO
that is claimed to exists; the set ${\cal R}_{a,w}$ is the set of its oracle rectangles.
We shall establish the claim  by induction on $k_a$.

If $k_a=1$, $a$ is a leaf. Take arbitrary rectangle $U_{a,w} \times V_{a,w}$. Either $i = {\sf lab}(a)$ is
correct on the rectangle, then $D_{a,w}$ is just the input $x_i$ and ${\cal R}_{a,w} = \emptyset$, or
not, and then $D_{a,w} = y_j$ and 
${\cal R}_{a,w} = \{R_j\}$
where $U_{a,w} \times V_{a,w}$ is $R_j$ in the enumeration (\ref{25.9.16a}).

Assume $k_a > 1$ and let $w \in \at$. For
$u \in U_{a,w}$ let $u^* \in \{0,1\}^{4^t}$ be a vector
whose bits $u^*_{\omega}$ are parameterized by
$\omega = (\omega_1, \omega_2) \in \at \times \at$
and such that: 
\begin{itemize}
\item $u^*_{\omega} = 1$ iff there is a $v \in V_{a,w}$
such that the communication of the players computing 
$S(u,v,a)$ evolves according to $\omega_1$ and the
computation of $S(u,v,a) \in_? F(u,v)$ evolves according to
$\omega_2$ (note that it has to end with the affirmation that 
$S(u,v,a) \in F(u,v)$.
\end{itemize}
Define $v^*_{\omega} \in \{0,1\}^{4^t}$ dually:
\begin{itemize}
\item $v^*_{\omega} = 0$ iff there is a $u \in U_{a,w}$
such that the communication of the players computing 
$S(u,v,a)$ evolves according to $\omega_1$ and the
computation of $S(u,v,a) \in_? F(u,v)$ evolves according to
$\omega_2$.
\end{itemize}
Let $U^*_{a,w}$ and $V^*_{a,w}$ be the sets of all these
vectors $u^*$ and $v^*$, respectively.

\medskip
\noindent
{\bf Claim 2:} {\em
There is a monotone formula $\varphi_{a,w}$ in $4^t$ variables $z_{\omega_1, \omega_2}$ 
and of size $2^{O(t)}$ separating
$U^*_{a,w}$ from $V^*_{a,w}$.}

\smallskip
\noindent

Claim 2 follows from a theorem of Karchmer and Wigderson \cite{KW}:
the players can find a coordinate $\omega$ in which $u^*_{\omega} =1$
and $v^*_{\omega} = 0$ by first computing $S(u,v,a)$
(getting thus $\omega_1$) and then deciding
$S(u,v,a) \in_? F(u,v)$ (obtaining thus $\omega_2$).
The strings $u \in U_{a,w}, v \in V_{a,w}$ yielding $u^*, v^*$ need not to be unique but that is not 
needed; it suffices that each player has a canonical way to pick one
such $u$ or $v$, respectively.

\medskip 

For $\omega_1 \in \at$ let $a_{\omega_1}$ be 
the node $S(u,v,a)$ computed
for some $u, v$ with communication  $\omega_1$. Then define 
a monotone circuit with a local oracle  by setting:
$$
D_{a,w} \ :=\ 
\varphi_{a,w}(\dots,\ 
z_{\omega_1, \omega_2}/D_{a_{\omega_1},\omega_2}
,\dots)\ 
$$
and:
$$
{\cal R}_{a,w} := \bigcup_{(\omega_1, \omega_2)} {\cal R}_{a_{\omega_1}, \omega_2}\ .
$$
As $k_{a_{\omega_1}} < k_a$, the induction hypothesis
implies that all $D_{a_{\omega_1},\omega_2}$
work correctly on all 
$U_{a_{\omega_1},\omega_2} \times V_{a_{\omega_1},\omega_2}$.
Thus,
by the definition of
the formula $\varphi_{a,w}$,
the circuit $D_{a,w}$ works also correctly.

\medskip

This concludes the proof of Claim 1 and of the theorem
(which follows from the claim by taking for $a$ the root of $G$).
The bound $\epsilon$ to the locality comes from our choice to start
with a protocol making an error for at most an $\epsilon$-part of $U \times V$.

\qed

It may be worthwhile to remark that the oracle rectangles of the CLO constructed in the
proof can be divided into $O(s)$ blocks (corresponding to different leaves)
such that the rectangles in each block are disjoint (they correspond to different
communication histories).

The particular CLO is constructed from a
particular $\rpp$ chosen by averaging. However, we could 
construct a CLO for each $\rpp$ and instead of estimating the locality
of the one CLO estimate the probability that a pair $(u,v)$ gets into an oracle
rectangle. We do not pursue this generality further here but we state
it formally as it may play a role in an eventual lower bound argument
for randomized protocols.

\begin{lemma} 
Assume that $\rrpp$ is a randomized protocol for $KW^m[U,V]$ of size $s$, communication
complexity $t$ and error $\epsilon$.

Then there is a distribution
$\rrcc$ over monotone circuits with a local oracle separating $U$ from $V$, each 
of size $s 2^{O(t)}$ and such that for any pair $(u,v) \in U \times V$:
$$
\prob_{{\bf r}}[(u,v)\ \mbox{ is in an oracle rectangle of } \rcc]\ \le \ \epsilon\ .
$$
\end{lemma}

The next two lemmas establish a form of converse of Lemma \ref{25.9.16b}.
Let $\umin$ be the set of $\le$-minimal elements of $U$ and $\vmax$ the set of
$\le$-maximal elements of $V$. In particular, no two elements of $\umin$ (or of $\vmax$),
respectively, are comparable and hence any partial Boolean function on $\umin$ (or on $\vmax$)
can be extended to a monotone one on $\umin \cup \vmax$.

\begin{lemma} \label{25.9.16c}
Assume $D, \{U_j \times V_j\}_{j \le e}$ is a monotone CLO separating $\umin$ from $\vmax$,
of size $s$ and locality $\mu$.

Then there is a protocol $\pro$ for $KW^m [\umin, \vmax]$ of size $s$, communication complexity
$2$ and making an error for at most $s \cdot \mu^{1/2}$-part of $\umin \times \vmax$. 
\end{lemma}

\prf

For each $j \le e$, the measure of $U_j \times V_j$ in $\umin \times \vmax$ is less than $\mu$
and hence

\begin{enumerate}

\item[(i)] either $|U_j|/|\umin| < \mu^{1/2}$,

\item[(ii)] or  $|V_j|/|\vmax| < \mu^{1/2}$.
\end{enumerate}
Define a monotone Boolean
function $f_j$ that is identically $1$ on $U_j$, identically $0$ on $V_j$, and for a string
from $\nn \setminus (U_j \cup V_j)$ it equals to $0$ in the case (i) or to $1$ in the case (ii).

Put $C(\overline x) := D(\overline x, f_1(\overline x), \dots, f_e(\overline x))$. Define a protocol
$\pro$ as follows:
\begin{itemize}

\item the vertices of $G$ are the nodes of $D$, the root is the output node and the edges
lead from a node of $D$ to its two input nodes,

\item for a node $a$ of $G$ corresponding to a subcircuit $E$ of $D$, define the consistency
condition by:
$$
a \in F(u,v)\ \mbox{ iff } \ 
(E(u, f_1(u), \dots, f_e(u)) = 1 \wedge E(v, f_1(v), \dots, f_e(v)) = 0)\ ,
$$

\item the strategy finds an input into $E$ that is also in $F(u,v)$,

\item the labeling $\sf lab$ assigns to input nodes $x_i$ of $D$ the value $i$
and to input nodes $y_j$ an arbitrary value, say $1$.

\end{itemize}
An error can occur only at the labeling of the input nodes corresponding to a
variable $y_j$. Because such a node is in $F(u,v)$, it must hold that $f_j(u) = 1$
and $f_j(v)=0$. 
In both cases (i) and (ii) considered in the definition of $f_j$
the measure of the rectangle of such pairs $(u,v)$ is
less than $\mu^{1/2}$ and there are at most $e \le s$ of them. This proves the lemma.

\qed

Next we show, for the sake of a completeness of the discussion, that one can
get a better estimate of the error of the protocol if one allows 
Boolean functions (and circuits) 
to have also a third value between $0$ and $1$. Denote the
third value $1/2$ and define the conjunction and the disjunction on $\{0, 1/2, 1\}$ as the
minimum and the maximum, respectively. Call such functions and circuits {\em 3-valued}.
We shall say that $D, \{U_j \times V_j\}_{j \le e}$ is a monotone 3-valued CLO separating
$U$ form $V$ if the condition 3. in the definition of the CLO is obeyed even 
w.r.t. to all monotone 3-valued functions $f_j$.

\begin{lemma} \label{25.11.16c}
Assume $D, \{U_j \times V_j\}_{j \le e}$ is a monotone 3-valued CLO separating $\umin$ from $\vmax$,
of size $s$ a locality $\mu$.

Then there is a protocol $\pro$ for $KW^m [\umin,\vmax]$ of size $s$, communication complexity
$2$ and making an error for at most $\mu$-part of $U \times V$. 
\end{lemma}

\prf

The construction of $\pro$ is similar to that in the proof of Lemma \ref{25.9.16c}
but we define the functions $f_j$ differently: $f_j$ equals to $1$ on $U_j$, to $0$ on
$V_j$ and to $1/2$ everywhere else.

With this definition the analysis at the end when an error occurs for a pair $(u,v)$ at 
a node corresponding to $y_j$ leads as before to a rectangle of $(u,v)$ such that
$f_j(u) = 1 \wedge f_j(v)=0$ but that is now simply $U_j \times V_j$. Hence the measure
of the set of pairs for which an error occurs is at most the locality of the CLO.

\qed

Let us conclude the section with a couple of remarks. 
The first one is that monotone CLOs simulate
efficiently monotone real circuits of \cite{Pud-cp} (circuits allowing any non-decreasing
real functions at gates); we shall show this 
in Lemma \ref{6.11.16a}. The second remark\footnote{I owe this
remark to Igor C. Oliveira.} is that general, non-monotone, CLOs are very strong:  
any two disjoint subsets of $\nn$ can be separated by a polynomial size CLO (in fact, a formula with
a local oracle) with
polynomially small locality. This is seen as follows: take the randomized protocol from Lemma \ref{24.11.16a}
and turn it into a non-monotone (dropping in the definition
the condition of monotonicity of oracle functions $f_j$)
CLO of size $poly(n, \epsilon^{-1})$ and locality $\epsilon$ separating $U$ from $V$ 
by the construction of Lemma \ref{25.9.16b}.

\section{Linear width} \label{4.9.15b}

The {\em linear width} of an $\rlin$-clause $C$ is the number of
$f$s in it; we shall denote it $\lw(C)$. For a set $\Phi$ of
$\rlin$-clauses denote by $\Phi \vdash_w C$ the fact that $C$ can be derived
in $\rlin$ from $\Phi$ by a proof whose all lines have linear width at most $w$.

When the linear width is small
the clauses have small
communication complexity (in the sense of  
Section \ref{20.3.16a}) and Theorem \ref{3.2} yields 
a small monotone protocol and that yields lower bounds
(cf. \cite[Sec.7]{Kra-interpol}).

Unfortunately, general $\rlin$ refutations need not to have small linear width.
It is easy to prove a lower bound on the linear width of an $\rlin$ refutation
by translating it into a polynomial calculus PC
refutation and by appealing to degree
lower bounds for that system. In particular, to an
$\rlin$-clause $C = \{f_1, \dots, f_k\}$ assign polynomial 
over $\ff$ $p_C := \Pi_{i \le k} (1-f_i)$: $C$ is satisfied by $a \in \nn$
iff $p_C(a)=0$.
An $\rlin$-refutation $\pi$ of a set $\Phi$ of $\rlin$ clauses
can be then straightforwardly translated into a PC refutation $\pi'$ of
the set of polynomials
$$
p_C\ ,\ \ C \in \Phi
$$
such that the degree of $\pi'$ is bounded above by the linear width of $\pi$.
In particular, the weakening rule and the binary rule translate into the
multiplication and the addition rules of PC, respectively.

To illustrate this lower bound argument let us consider as a specific
example the set $\neg PHP_n$ of $\rlin$ clauses:

\begin{itemize}

\item $1 - x_{i j}, x_{k j}$, for $i \neq k$ and any $j$,

\item $1 - x_{i j}, x_{i k}$, for any $i$ and $j \neq k$,

\item $\sum_j x_{i j}$, any $i$,

\end{itemize}
with variables $x_{i j}$, $i \in [n+1], j\in [n]$. The linear width of these clauses is $1$.
However, the set of polynomials $p_C$ for $C \in \neg PHP_n$ is precisely the set
for which the degree $n/2$ lower bound for PC refutations was established by Razborov
\cite{Raz-pc}.

We shall employ the approximation method in Section \ref{4.9.15c}
to reduce in a sense the linear width. 
This construction introduces, however, some error into 
derivations (modelled in one of the constructions by new initial clauses
to be called $\newax$) and this prevents the simple reduction to 
PC described above.

\section{Randomized feasible interpolation for\\
 $\rlin$} \label{4.9.15c}

In this section we use the Razborov-Smolensky approximation method \cite{Raz87,Smo}
to reduce in a sense the linear width of not too large $\rlin$ refutations.

\begin{theorem} \label{30.10.16e}
Assume the set-up conditions (1)-(7)
and assume that sets\newline
 $A_1, \dots, A_m, B_1, \dots, B_\ell$ are defined by $\rlin$-clauses. 

Let $\pi$ be an $\rlin$ refutation of
(the clauses defining) these sets with $k$ steps. 
Let $w \geq 1$ be any parameter.

Then  there is a randomized protocol for $KW^m[U,V]$) of size at most $k+n$,
communication complexity $O(w \log n)$ and of error
at most $3 \cdot 2^{-w} k$.

Moreover, if the refutation $\pi$ is tree-like
then also $G$ is tree-like.
\end{theorem}

\prf 

Let $D$ be any $R(LIN)$-clause, i.e. a clause formed by some linear polynomails. 
Following \cite{Raz87,Smo}
define a $2^{-w}$-approximation $(Y_{\bs})_{\bs}$ of $D$ by the following process:
\begin{itemize}

\item Using the sample $\bs$ pick independently at random $L_1, \dots, L_{w} \subseteq D$,

\item put $Y_{\bs}$ to be the set defined by $\bigvee_{j \le w} \sum L_j$,

\end{itemize}
($\sum L_j$ is the sum of all linear polynomials in $L_j$). 

\medskip
\noindent
{\bf Claim:} {\em Let $D$ be an 
$\rlin$-clause of linear width $w$. Then $MCC_U(D) = O(w \log n)$.}

\medskip

Let us write the $w$ linear functions forming $D$ in a matrix form as:
$$
A x + B y + C z + E\ .
$$
The U-player sends $A u$ and $B q^u$ and the V-player sends $A v$ and $C r^v$, $4w$ bits in total.
After this they know the truth values of $(u, q^u, r^v) \in D$ and 
$(v, q^u, r^v) \in D$ and if they differ they can use the binary search on a differing
row in $A u$ and $A v$ to find $i$ for which $u_i \neq v_i$ ($2 log n$ bits in total).

It remains to estimate the communication complexity of the
task 4. from the definition of $MCC_U$ under the assumption
that  $(u, q^u, r^v) \in D$ and $(v, q^u, r^v) \notin D$, i.e.:
$$
A u + B q^u + C r^v + E \neq \overline 0\ \mbox{ and }\ 
A v + B q^u + C r^v + E = \overline 0\ .
$$
In particular, $A u \neq Av$.

The players will attempt to put $A$ in a reduced-row echelon form but by a specific process. 
The U-player sends $i_1 \in [n]$ ($\log n$ bits) such that $u_{i_1}=0$ and the $x_{i_1}$-column
in $A$ is non-zero. 
The players then both separately transform 
$A$ using the elementary row and column operations in some canonical 
way to a unique matrix $A^1$ whose first column corresponds to $x_{i_1}$ and $A^1_{1,1} = 1$ and all
other entries in the first column are $0$.

In the second step they apply the same process to $A^1$, not using $x_{i_1}$. That is, 
the U-player sends $\log n$ bits identifying some
$i_2 \in [n]$, $i_2 \neq i_1$, such that $u_{i_2} = 0$ and 
the $x_{i_2}$-column in $A^1$ has a non-zero element in one of the 
rows $2, \dots, w$. Then they again separately transform $A^1$ into 
$A^2$ with the first two columns corresponding to $x_{i_1}$ and $x_{i_2}$ and the left-upper
corner $2 \times 2$ submatrix being the identity matrix $I_2$ and all other entries in the first two
columns being $0$.

They proceed analogously as long as it is possible. Two cases may occur:
\begin{enumerate}

\item [(i)] After $t \le w$ steps $A^t$ is in the row-reduced echelon form: the 
left-upper
corner $t \times t$ submatrix being the identity matrix $I_t$ and all other entries in the first $t$
columns being $0$, and all rows $t+1, \dots, w$ are zero.

\item [(ii)] After some step $t < w$ $A^t$ is not in the row-reduced echelon form
but the U-player has nothing to choose: there is no 
$i \neq i_1, \dots, i_t$ such that the $x_i$-column in $A^t $ has a non-zero element
in one of the rows $t+1, \dots, w$ and $u_{i} = 0$.

\end{enumerate}
In Case (i) we can switch the values of some $u_i$, $i \in \{i_1, \dots, i_t\}$,
from $0$ to $1$ to get $u'\geq u$ such that $A^t u'= A^t v$ and hence
$(u', q^u, r^v) \notin D$. 

In Case (ii) the rows $t+1, \dots, w$ need not to be zero but
$A_{i j} \neq 0$ for $i, j > t$ implies that $u_i = 1$ (thinking of 
the $i$-th column as corresponding to $x_i$). If for one such $i$ $v_i =0 $, the V-player sends
the $\log n$ bits to identify it; they found $i$ such that $u_i = 1 \wedge v_i = 0$.
If all such $v_i = 1$ then $G u = G v$ where $G$ is the $(w-t)\times n$ matrix consisting of the
last $w-t$ rows of $A^t$.
Writing the first $t$ rows of $A^t$ as $(I_t, H)$, where $H$ is a $t \times (n-t)$ matrix, we
see we can find some $u'\geq u$ changing only some $u_i$, $i \in \{i_1, \dots, i_t\}$,
from $0$ to $1$
such that $(I_t,H) u' = (I_t, H) v$ and hence also $A^t u'= A^t v$ and 
$(u', q^u, r^v) \notin D$. 

In all cases the players solved the task 4. and they exchanged $O(w \log n)$ bits at most.

\medskip

Applying Theorem \ref{29.3.16b}  concludes the proof of the theorem.

\qed

We now give an alternative proof of the randomized feasible interpolation for $\rlin$,
referring to Theorem \ref{29.3.16a} this time. It is more laborious and gives somewhat worse bounds
on the size of the resulting protocols but 
it may be useful in connections with the problem of resolution over low degree
polynomial calculus that we shall discuss in the Section \ref{4.9.15e}, and it 
also puts 
$\rlin$ in a direct relation with the random R of \cite{BKT} 
(see Section \ref{4.9.15e}).

Let $\pi$ be an $\rlin$ refutation
of $\Phi := A_1, \dots, A_m, B_1, \dots, B_{\ell}$ and let
$w \geq 1$ be a parameter to be specified later. In this
situation we perform the following random process ${\bf r}$ and transform
$\pi$ to an $\rlin$ refutation $\pi({\bf r})$ of $\Phi$ extended by
a set $Ax(\pi, {\bf r})$ of extra clauses:

\begin{enumerate}

\item For each $C \in \pi$ pick independently at random subsets 
$L_1, \dots, L_\ell \subseteq C$ and
form clause $C^{\bf r} := \{\sum L_1, \dots, \sum L_\ell \}$.

\item For each $C \in \pi$, $C = f_1, \dots, f_k$, add to the set $Ax(\pi, {\bf r})$ 
the following $k$ clauses:
$$
C^{\bf r}, f_j+1\ \mbox{ , for } j = 1, \dots, k\ .
$$

\item Transform $\pi$ into $\pi({\bf r})$, following the construction below, summarized
in Lemma \ref{4.4.15a}.

\end{enumerate}
Clauses in 2. formalize that $f_j=1$ implies that $C^{\bf r} = 1$.
Before we describe $\pi({\bf r})$ we need to establish a few simple
facts.

\medskip
\noindent
{\bf Claim 1:} {\em 
For all assignments $a \in \nn$: $C^{\bf r}(a)=1$ implies
$C(a)=1$. For any $a \in \nn$ the probability that $C^{\bf r}(a)=0 \wedge C(a)=1$
is at most $2^{-w}$.}

\medskip
\noindent
{\bf Claim 2:} {\em (a) For any $g,h$: $g+h \vdash_2 g,h$.

(b) For any $C \in \pi$ and $g \in C^{\bf r}$:
$g \vdash_{|C|} C$.}

\smallskip
In part (a): derive from $g+h$ clause $g,g+h$ and also an $\ff$-axiom $g, g+1$
from which $g,h$ follows by the binary rule and contraction. In part (b): if 
$g = f_{j_1} + \dots + f_{j_v}$ use part (a) to derive from $g$
clause $f_{j_1} + \dots + f_{j_{v-1}}, 
f_{j_v}$, and then repeat this to remove from the sum all $f_j$s to
get the clause $f_{j_1}, \dots, f_{j_v}$ from which $C$ follows by the
weakening rule. 

\medskip
\noindent
{\bf Claim 3:} {\em
Let $C \in \pi$, $C = f_1, \dots, f_k$, and let $g = f_{j_1} + \dots + f_{j_v}$
be an arbitrary sum of a non-empty subset of $C$ (i.e. not
necessarily in $C^{\bf r}$). Then
$$
Ax(\pi, {\bf r}), \{g\} \ \vdash_{w + 3} C^{\bf r}\ . 
$$
}

\smallskip

By Claim 2(a) derive in linear width $2$ from $g$ clause 
$f_{j_1} + \dots + f_{j_{v-1}}, f_{j_v}$ and combine this by the binary
rule and contraction with clause 
$f_{j_v}+1, C^{\bf r}$ from $\newax$ to get
$$
f_{j_1} + \dots + f_{j_{v-1}}, C^{\bf r}
$$
in linear width bounded by $w + 3$. Then repeat the same process
to remove from the sum polynomials $f_{j_{v-1}}, f_{j_{v-2}}, \dots,
f_{j_1}$ to end up just with $C^{\bf r}$.

\medskip
\noindent
{\bf Claim 4:} {\em
Assume 
$$
\frac{C}{C,h}
$$
is an inference in $\pi$.
Then 
$$
\newax, C^{\bf r}\ \vdash_{2 w +2}\ (C,h)^{\bf r}\ .
$$
}

\smallskip
Assume $C^{\bf r} = \{g_1, \dots, g_{w}\}$ where each $g_i$ is
a sum of some polynomials from $C$ and thus also from $C,h$. So repeating Claim 3
$w$-times to remove $g_{w}, g_{w -1}, \dots, g_1$ we derive
$(C,h)^{\bf r}$. The linear width is at most $w+3$ (from Claim 3) plus 
$w - 1$ (for side polynomials $g_1, \dots, g_{w-1}$), i.e. at most
$2w + 2$ in total.

\medskip
\noindent
{\bf Claim 5:} {\em
Assume
$$
\frac{C,g\ \ \ \ C,h}{C,g+h+1}
$$
is an inference in $\pi$. Then 
$$
\newax, (C,g)^{\bf r}, (C,h)^{\bf r})\ \vdash_{2w + 3}\ (C,g+h+1)^{\bf r}\ .
$$
}

\smallskip

We proceed as in Claim 4 and attempt to derive from $\newax, (C,g)^{\bf r}$ clause
$(C,g+h+1)^{\bf r}$. The only obstacle to doing so is when the polynomial $g$
occurs in a sum in $(C,g)^{\bf r}$: in that case we leave it as a side polynomial. 
That is, from $(C,g)^{\bf r}$ we derive $(C,g+h+1)^{\bf r}, g$ in linear width
at most $2w+2+1 = 2w + 3$.

Analogously from $(C,h)^{\bf r}$ derive 
$(C,g+h+1)^{\bf r}, h$ and then by the binary rule 
$$
(C,g+h+1)^{\bf r}, g+h+1\ .
$$
From that we get the wanted $(C,g+h+1)^{\bf r}$ using the axiom
$$
(g+h+1)+1, (C,g+h+1)^{\bf r}
$$
from $\newax$, the binary rule and a contraction.

\bigskip

The following lemma follows form the last two claims.

\begin{lemma} \label{4.4.15a}
Let $\pi$ be an $\rlin$ refutation of $A_1, \dots, A_m, B_1, \dots, B_\ell$ 
consisting of $k$ clauses
and of linear width $w_0$. 
Let $w \geq 1$ be an arbitrary parameter.
Then for a random
$\bf r$ there is an $\rlin$-refutation $\pr$ of 
$$
\Phi, \newax 
$$
of linear width bounded above by
$$
w'\ :=\ 2w +3
$$
and with at most $O(w w_0 k)$ clauses.

\end{lemma}

\prf

The bound to the linear width follows from the last two claims, using also that
$$
\Phi, \newax\ \vdash_{w'} \Phi^{\bf r}\ .
$$
The bound to the number of clauses follows by inspecting that in both
Claims 4 and 5 the constructed derivations have $O(w w_0)$ clauses.

\qed

We used in this construction the syntactic version of $\rlin$ rather than the semantic one
in order to generate explicitly the sets $\newax$.

Now we can apply Theorem \ref{29.3.16a}. The values of parameters appearing
in that theorem are:
\begin{itemize}

\item $\epsilon := 2^{-w} k$: the conjunction of
axioms in $\newax$ corresponding to any one clause
in $\pi$ are $2^{-w}$-approximations of $\NN$ (Claim 1).

\item Number of steps: $O(w w_0 k)$.

\item Monotone communication complexity: $O(w \log n)$. 

\end{itemize}

\begin{theorem} \label{30.10.16d}
Assume the set-up conditions (1)-(7)
and assume that sets \newline
$A_1, \dots, A_m, B_1, \dots, B_\ell$ are defined by $\rlin$-clauses.

Let $\pi$ be an $\rlin$-refutation of (the clauses defining)
these sets with $k$ steps and of the linear width bounded
by $w_0$.

Then for every $w \geq 1$
there is a randomized protocol $\rrpp$ for $KW^m[U, V]$ of size at most $O(w w_0 k) + n$,
communication complexity $O(w \log n)$ and of error
at most $2^{-w}k$.

Moreover, if the refutation $\pi$ is tree-like
then also $G$ is tree-like.

\end{theorem}

Using Lemma \ref{25.9.16b} we can turn Theorems \ref{30.10.16e}
and \ref{30.10.16d}
into statements about separating monotone
CLOs (we use Theorem \ref{30.10.16e} in the corollary).

\begin{corollary} \label{30.10.16f}
{\ }

Assume the set-up conditions (1)-(7)
and assume that sets\newline
 $A_1, \dots, A_m, B_1, \dots, B_\ell$ are defined by $\rlin$-clauses. 

Let $\pi$ be an $\rlin$ refutation of
(the clauses defining) these sets with $k$ steps. 
Let $w \geq 1$ be any parameter.

Then  there is a monotone CLO of size at most $(k+n) 2^{O(w \log n)}$ and of locality
at most $3 \cdot 2^{-w} k$ separating $U$ from $V$.

Moreover, if the refutation $\pi$ is tree-like
then the monotone CLO is a formula.

\end{corollary}

\section{Randomized feasible interpolation for CP}   \label{20.3.16c}

Following \cite{Kra-game} call a
semantic derivation {\em CP-like} iff the 
proof steps are defined by integer linear inequalities.
CP-like derivations were interpolated in \cite{Kra-game} by
protocols but their complexity was measured in terms of the real game
defined there: players send each a real number to a referee and he
announces how are these ordered. 
The real communication complexity of a multi-function $R$, $\ccr(R)$, is the minimal number of
rounds (of sending numbers to the referee in an optimal protocol) needed to compute a valid
value for $R$ in the worst case.
We can use this notion to measure the communication complexity of our protocols $\mathbf P$ and 
define $\ccr({\mathbf P})$ analogously to how $\cc({\mathbf P})$ was defined.  
We will not recall details 
as we will use here only the relation of the
real communication complexity to the well-established probabilistic communication complexity.

Let $R$ be a multi-function defined on $U \times V \subseteq \nn \times \nn$ 
and let $\pcc(R)$ be the 
probabilistic communication complexity of
a multi-function $R$ with public coins
and error $\epsilon > 0$. The following equality was derived in
\cite[L.1.6]{Kra-game} from a result of Nisan \cite{Nis93}.
For $\epsilon < \frac{1}{2}$ it holds
\begin{equation}\label{29.3.16c}
\pcc (R) \le \ccr (R) \cdot O(\log n + \log \epsilon^{-1})\ .
\end{equation}
We will use \cite[Thm.3.3]{Kra-game}.

\begin{theorem} \label{6.11.16b}
Assume the set-up conditions (1)-(7). Assume that the
sets $A_1, \dots, A_m$ and $B_1, \dots, B_\ell$ are defined by
integer linear inequalities and that there is a CP-like refutation
$\pi$ of $\tilde A_1, \dots, \tilde A_m, \tilde B_1, \dots, \tilde B_{\ell}$
that has $k$ steps.

Then for any $\epsilon < 1/k$ there is a randomized protocol for $KW^m[U,V]$ of size 
$k+n$, communication complexity 
$O(\log (n/\epsilon))$
and of error at most $\epsilon k$.

Moreover, if the refutation $\pi$ is tree-like
then also $G$ is tree-like.
\end{theorem}

\prf

Theorem 3.3. of \cite{Kra-game} shows that there is a protocol
for $KW[U,V]$ (resp. for $KW^m[U,V]$) of the stated size and with the real
communication complexity $O(1)$. Then (\ref{29.3.16c}) implies
that that protocol can be simulated by a randomized protocol 
of communication complexity $O(\log (n/\epsilon))$
which, for given $u,v$,
computes at every node $x$ the strategy function and the consistency condition with
error at most $\epsilon$. 
Hence the total error is estimated by $\epsilon k$.
This entails the theorem.

\qed

Note that analogously to Corollary \ref{30.10.16f}
this can be turned into a statement about separating monotone CLOs.
However, it is more direct to use the argument from the preceding proof
to show that monotone CLOs efficiently simulate monotone 
real circuits of Pudl\' ak \cite{Pud-cp} which do separate pairs $U, V$ by the interpolation theorem established there.

\begin{lemma} \label{6.11.16a}
Assume $U, V \subseteq \nn$ and $U$ is closed upwards (or $V$ downwards). 
Let $C$ be a monotone real circuit of size $s$ separating $U$ from $V$.

Then for every $0 < \epsilon < \frac{1}{2}$ there is a monotone CLO $D$
separating $U$ from $V$, having size $s (\frac{n}{\epsilon})^{O(1)}$
and locality $\mu \le s \epsilon$.

In particular, for any $\mu > 0$ there is a  
monotone CLO separating $U$ from $V$ with locality $\le \mu$ and size 
$(n s \mu^{-1})^{O(1)}$.
\end{lemma}

\prf

Circuit $C$ yields a protocol for $KW^m[U,V]$ of size $s$ and real communication
complexity $O(1)$: the graph of the protocol is $C$ turned upside down (output is the root), 
the consistency
condition $F(u,v)$ consists of subcircuits $E$ where $E(u) > E(v)$, and the strategy
is defined so that the consistency condition is preserved.

 As in the proof of Theorem \ref{6.11.16b} the protocol can be turned into 
a randomized protocol of size $s$, communication complexity $O(\log(n/\epsilon))$ 
and error at most $s \epsilon$. The required monotone CLO then exists by Lemma 
\ref{25.9.16b}.

The particular case is obtained by setting $\epsilon := s/\mu$.

\qed

Let us remark that the constructions underlying Theorem \ref{6.11.16b} and Lemma 
\ref{6.11.16a} apply also to the proof system R(CP) of \cite{Kra-modules}
operating with clauses formed by CP-inequalities and yield a small separating CLO for
small width. In particular, if each clause in an R(CP)-refutation 
has size at most $w$ then the (monotone) real communication complexity is
at most $w$ and this yields a monotone separating CLO of the size as in Lemma 
\ref{6.11.16a} for 
$w = O(\log (n/\epsilon))$.

\section{The lower bound problem for monotone CLOs} \label{20.3.16d}

This section is devoted to a discussion of the problem to establish
a lower bound for monotone circuits with
a local oracle  separating two sets $U$ and $V$ (obeying all set-up conditions
(1) - (7)). This would imply
via Lemma \ref{25.9.16b} also a lower bound for randomized protocols
for $KW^m[U,V]$ and hence a length-of-proofs lower bound for $\rlin$.

We shall consider the classical pair of disjoint sets of 
graphs having a large clique and of graphs colorable by a small number of colors.
Let $n_0 \geq \omega > \xi \geq 1$ and put $n := {n_0 \choose 2}$.
We shall
identify in this context $[n]$ with the set of unordered pairs of distinct elements from
$[n_0]$; we think of each such pair as denoting a potential edge in
a graph with vertices $[n_0]$.

Take for $U \subseteq \nn$ the set $Clique_{n_0,\omega}$ of all graphs on
$[n_0]$ that contain a clique of size $\omega$. We shall also denote
by $Clique_{n_0, \omega}(p,q)$ the set of the following clauses in
atoms $p_{i j}$, $i\neq j \in [n_0]$, and
$q_{u i}$, $u = 1, \dots, \omega$ and $i \in [n_0]$ (hence there are 
$s = \omega \cdot n_0$ $q$-atoms):
\begin{itemize}
\item $\bigvee_{i \in [n_0]} q_{u i}$, one for each $u \in [\omega]$,
\item $\neg q_{u i} \vee \neg q_{v i}$, one for all $u < v \in [\omega]$
and $i \in [n_0]$,
\item $\neg q_{u i} \vee \neg q_{v j} \vee p_{i j}$, 
one for all $u < v \in [\omega]$ and $i \neq j \in [n_0]$.
\end{itemize}
Sets $A_i$ from the set-up condition (2)
are the sets defined by these clauses.

The set $V \subseteq \nn$ will be the set of graphs on $[n_0]$
that are $\xi$-colorable. We shall
denote it $Color_{n_0,\xi}$ and by $Color_{n_0,\xi}(p, r)$ the set
of the following clauses in the $p$-atoms and atoms
$r_{i a}$, $i\in [n_0]$ and $a \in [\xi]$ (there are $n_0 \cdot \xi$ $r$-atoms):
\begin{itemize}
\item $\bigvee_{a \in [\xi]} r_{i a}$, one for each $i \in [n_0]$,

\item $\neg r_{i a} \vee \neg r_{i b}$, one for all $a < b \in [\xi]$
and $i \in [n_0]$,

\item $\neg r_{i a} \vee \neg r_{j a} \vee \neg p_{i j}$, 
one for all $a \in [\xi]$ and $i\neq j \in [n_0]$.

\end{itemize}
Sets $B_j$ from the set-up condition (2)
are the sets defined by these clauses.

If we identify a truth assignment $w \in \nn$ to the $p$-atoms with graph $G_w$ on $[n_0]$,
truth assignments to $q_{u i}$ satisfying
$Clique_{n_0, \omega}(w, q)$ correspond to injective (multi-)maps from $[\omega]$
onto a clique in $G_w$ and analogously truth assignments
to $r_{i a}$ making $Color_{n_0, \xi}(w, r)$ true
correspond to colorings of $G_w$ by $\xi$
colors. Thus if $ \omega > \xi $ the sets $U$ and $V$ are disjoint
and its is easy to see that they, together with the clauses above,
satisfy the set-up conditions (1)-(7) from Section \ref{20.3.16a}.

\bigskip

Let us first note that a lower bound for a monotone CLO with oracle rectangles inside
$U \times V$
can be derived as an easy consequence of a theorem of Jukna \cite[Thm.3]{Juk-local},
generalizing an earlier result by Yao \cite{Yao-local}.
In particular, \cite[Thm.3]{Juk-local} states that there is no small (polynomial size)
monotone circuit computing the characteristic function $\chi_U$ of $U$
for $\omega = (n_0/\log n_0)^{2/3}$ even if the circuits are allowed to use
at gates arbitrary monotone Boolean functions as long as all their min-terms have size
$o(\omega)$. In the case of a monotone CLO with oracle rectangles $U_j \times V_j$
we can take for all functions $f_j$ the disjunction $f$ of all conjunctions
\begin{equation}\label{28.9.16a}
\lceil X \rceil \ :=\ \bigwedge_{i \neq j \in X} p_{i j}
\end{equation}
where sets $X \subseteq [n_0]$ run over all sets of vertices of size $\xi + 1$.
Clearly $f$ is identically $1$ on $U$ and $0$ on $V$ and hence
if, say, $\xi = \omega^{1/2}$, Jukna's \cite[Thm.3]{Juk-local} applies. 
However, this is not good enough: we want a stronger lower bound but more importantly
we need a lower bound for monotone CLOs separating $U$ from $V$ and not just for those
computing $\chi_U$.

The classical result of Alon and Boppana \cite{AloBop}, strengthening Razborov's 
\cite{Raz85} lower bound, offers such a lower bound for ordinary monotone circuits.
 
\begin{theorem} [Alon and Boppana {\cite[Thm.3.11]{AloBop}}] 
{\ }

Assume that $3 \le \xi < \omega$ and $\sqrt{\xi} \omega \le
\frac{n_0}{8 \log n_0}$ . Then any monotone circuit separating 
$Clique_{n_0,\omega}$ from $Color_{n_0, \xi}$
must have the size at least 
$2^{\Omega(\sqrt{\xi})}\ .$
\end{theorem}

It appears possible that the same lower bound holds also for monotone CLOs 
with a small constant locality. 
Alluding to Boppana and Sipser \cite[L.4.2]{BopSip} we prove at least the following
partial result for monotone CLOs of the restricted form
\begin{equation} \label{25.11.16a}
D\ :=\ \bigvee_{i \le a} (\lceil X_i \rceil \wedge C_i(\overline y))
\end{equation}
where 
\begin{enumerate}

\item[12.1] $|X_i| \le \lfloor \xi^{1/2}\rfloor$ and $\lceil X_i \rceil$ is defined as in 
(\ref{28.9.16a}) using variables $x_{i j}$ in place of $p_{i j}$,

\item[12.2] $C_i(\overline y)$ is a monotone circuit of an arbitrary size not containing the $x$-variables,

\item[12.2] the size $a$ of the disjunction is arbitrary.
\end{enumerate}

\begin{lemma}  \label{3.10.16b}

Assume that $4 \le \xi < \omega$ and that $n_0$ is large enough. 
Then no monotone circuit with a local oracle $D$ of the form (\ref{25.11.16a}),
satisfying conditions 12.1-3 
and with locality $\mu \le \frac{1}{16}$
separates 
$Clique_{n_0,\omega}$ from $Color_{n_0, \xi}$.
\end{lemma}

The proof of the lemma will be summarized after Lemma \ref{25.11.16b}. 

\bigskip

A CLO separating $U (=\clique)$ from $V (=\color)$ separates also $\umin$ from $\vmax$. Note that
elements of $\umin$ are graphs consisting of a clique of size $\omega$ and having no other edges
and elements of $\vmax$ are $\xi$-partite graphs with all possible edges among the different parts. 
These two sets are
called in \cite{AloBop,BopSip} {\em positive} and {\em negative} examples, respectively.
In fact, for the counting purposes the negative examples are represented as $\xi$-colorings of
$[n_0]$, each coloring determining the maximal graph for which it is still a graph coloring.

Let $D(\overline x, \overline y), {\cal R}$ be a monotone CLO 
of the form (\ref{25.11.16a}), satisfying 12.1-3, with locality $\mu$ and
with $e$ oracle rectangles $U_j \times V_j$.
Let 
$$
Bad\ :=\ \bigcup_{j \le e} U_j \times V_j\ \subseteq \ \umin \times \vmax .
$$
We know that $|Bad| \le \mu \cdot |\umin \times \vmax|$.

In the argument 
we shall consider other rectangles inside $\umin \times \vmax$ and $y$-variables attached to them.
Let us introduce the following notation. For $U'\subseteq \umin$ and $V' \subseteq \vmax$
let $y[U',V']$ be a new variable. Its {\em valid interpretation} is any monotone Boolean function
$h : \nn \rightarrow \bits$ that is $1$ on $U'$ and $0$ on $V'$. Two specific valid interpretations
of the $y$-variables are:

\begin{itemize}

\item ${\cal F}_U$-interpretation: each $y[U',V']$ is interpreted by the Boolean function
that is $1$ on $U'$ and $0$ everywhere else on $\umin \cup \vmax$,

\item ${\cal F}_V$-interpretation: each $y[U',V']$ is interpreted by the Boolean function
that is $0$ on $V'$ and $1$ everywhere else on $\umin \cup \vmax$,

\end{itemize}
(we only care for values on $\umin \cup \vmax$).
Let $E(\overline x, \overline y)$ be a monotone circuit involving also some of the $y$-variables
and let ${\cal F}$ be a valid interpretation of the $y$-variables. Then 
$$
E(\overline x, {\cal F})
$$
denotes the Boolean function obtained by substituting for each $y$-variable in $E$ the function
interpreting it in $\cal F$.

\begin{lemma} \label{9.11.16a}
Let $E(\overline x, \overline y)$ be a monotone circuit. It holds on $\umin \cup \vmax$:

\begin{enumerate}

\item For any valid interpretation $\cal F$:
$$
E(\overline x, {\cal F}_U) \le E(\overline x, {\cal F}) \le E(\overline x, {\cal F}_V)\ .
$$

\item For ${\cal F} = {\cal F}_U, {\cal F}_V$:
$$
(y[U_1,V_1] \vee y[U_2,V_2])({\cal F})\ =\ 
y[U_1 \cup U_2, V_1 \cap V_2]({\cal F})\ .
$$

\item For ${\cal F} = {\cal F}_U, {\cal F}_V$:
$$
(y[U_1,V_1] \wedge y[U_2,V_2])({\cal F})\ =\ 
y[U_1 \cap U_2, V_1 \cup V_2]({\cal F})\ .
$$

\item If both $U_1 \times V_1$ and $U_2 \times V_2$ are subsets of $Bad$, so are
$U_1 \cup U_2 \times V_1 \cap V_2$ and $U_1 \cap U_2 \times V_1 \cup V_2$.
\end{enumerate}

\end{lemma}

\prf

Parts 1 and 4 are obvious. Let $\chi_W$ be the characteristic function of $W \subset \nn$.
For Part 2:
$$
(y[U_1,V_1] \vee y[U_2,V_2])({\cal F}_U) = \chi_{U_1} \vee \chi_{U_2} =
\chi_{U_1 \cup U_2} = y[U_1 \cup U_2, V_1 \cap V_2]({\cal F}_U)
$$
and
$$
(y[U_1,V_1] \vee y[U_2,V_2])({\cal F}_V) = \chi_{\setminus V_1} \vee \chi_{\setminus V_2} =
\chi_{\setminus(V_1 \cap V_2)} = y[U_1 \cup U_2, V_1 \cap V_2]({\cal F}_V)\ .
$$
Part 3 is analogous.

\qed

We shall argue that
either $D(\overline x, {\cal F}_U)$ rejects a lot of $\umin$ or that
$D(\overline x, {\cal F}_V)$ accepts a lot of $\vmax$. 
The choice to evaluate
how well $D$ works on $\umin$ using the interpretation ${\cal F}_U$
and on $\vmax$ using ${\cal F}_V$ gives us (due to Part 1 of Lemma \ref{9.11.16a})
the best chance to detect errors.

Note that $\lceil X \rceil$ is equivalent to
$$
\lceil X \rceil y[\umin, \emptyset]
$$
under the two extreme interpretations as
$y[\umin, \emptyset]$ is $1$ on $\umin$
under ${\cal F}_U$ and $1$ on
both $\umin$ and $\vmax$ under ${\cal F}_V$. So we could have allowed in (\ref{25.11.16a}) also
stand-alone terms $\lceil X \rceil$ and if we defined $\lceil \emptyset \rceil := 1$ also
stand-alone $y$-variables.

\begin{lemma} \label{25.11.16b}
Assume $\mu \le 1/16$. Then for any monotone CLO $E$ of the form 
$$
E = \bigvee_{i \le a} (\lceil X_i \rceil \wedge y[U_i, V_i])
$$
where $a$ is arbitrary, $|X_i| \le \lfloor \xi^{1/2} \rfloor$ and 
all rectangles $U_i \times V_i$ are  subsets of $Bad$ it holds:

\begin{enumerate}

\item  Either $E({\cal F}_V)$
accepts at least $1/4$ of $\vmax$,

\item or $E({\cal F}_U)$ rejects at least $3/4$ of $\umin$.

\end{enumerate} 

\end{lemma}

\prf

If $E$ is the empty disjunction, it is constantly zero and the second option
occurs. 

If not, note that as all rectangles $U_i \times V_i $ are subsets of $Bad$, their measure 
in $\umin \times \vmax$
at most $\mu$. Hence at least one of its sides $U_i$ or $V_i$
has the measure at most $\mu^{1/2}$ in $\umin$ or $\vmax$, respectively. Now
consider two cases:

\begin{enumerate}

\item There is a term $\lceil X_i \rceil \wedge y[U_i, V_i]$ in $E$
with $V_i$ having
the measure at most $\mu^{1/2}$ in $\vmax$,

\item not 1.
\end{enumerate}
Denote $\ell := \max_{j \le a} |X_j|$; we have $\ell \le \lfloor \xi^{1/2}\rfloor$.

In the first case the term 
$\lceil X_i \rceil \wedge y[U_i, V_i]({\cal F}_V)$ accepts at least the fraction of
$$
[1 - \frac{{\ell\choose 2}}{\xi}] - \mu^{1/2} \geq 
[\frac{3}{4} - \frac{{\ell\choose 2}}{\xi}] \geq \frac{1}{4}
$$
elements $v \in \vmax$: the first term is the same estimate as in \cite[L.4.2]{BopSip}, 
the second accounts for the elements of $V_i$.

In the second case use ${\cal F}_U$: 
all $y[U_i,V_i]({\cal F}_U)$ are $1$ only inside $U_i$ and hence $E$
accepts at most the subset  $\bigcup_i U_i$ of $\umin$. But for each $u$ from this union 
the pair $(u,v) \in Bad$ for at least a fraction of $\mu^{1/2}$ of elements $v$ of $\vmax$.
Hence the measure of the union
is at most $\mu^{1/2} \le \frac{1}{4}$.

\qed

Now we can derive Lemma \ref{3.10.16b}. By parts 2 and 3 of Lemma \ref{9.11.16a},
each subcircuit $C_i(\overline y)$ of $D$ is equivalent under both ${\cal F}_U$ and ${\cal F}_V$
to some $y[U_i, V_i]$ such that, by part 4 of that lemma, $U_i \times V_i \subseteq Bad$.
Hence Lemma \ref{25.11.16b} applies.

\bigskip

Let us remark that there is a certain discrepancy in the sizes when protocols are turned to CLOs
in Lemma \ref{25.9.16b} and CLOs are transformed into protocols in Lemma \ref{25.9.16c}.
Thus even if the lower bound for monotone CLOs was not valid one could still try the tight
$3$-valued version of Lemma \ref{25.11.16c}.

\section{Concluding remarks} \label{4.9.15e}

We remark without elaborating it that  Theorem \ref{29.3.16a} yields a randomized feasible
interpolation\footnote{A different one than \cite{Kra-rres}.} for  
the random resolution system proposed informally by Dantchev and defined formally
by Buss, Kolodziejczyk and Thapen \cite[Sec.5.2]{BKT}.
Pudl\' ak and Thapen \cite{PudTha} 
consider more variants of the definition and they prove a feasible interpolation
for the tree-like case. According to the definition from Buss et.al. \cite{BKT}
an $\epsilon$-random resolution refutation distribution
of a set of clauses $\Phi$ is a random distribution $(\pi_{\rr})_{\rr}$
of resolution refutations of $\Psi \cup \Delta_{\rr}$, where
$\Delta_{\rr}$ are sets of clauses such that any fixed
truth assignment fails to satisfy  $\bigwedge \Delta_{\rr}$ with the probability at most
$\epsilon$.
In other words, if $X_{\rr}$ is the set of assignments satisfying
all clauses in $\Delta_\rr$ then $(X_\rr)_\rr$ is an $\epsilon$-approximation
of the universe of all assignments. 
The number of steps in such a random refutation is the maximal number of
steps among all $\pi_{\rr}$.

\medskip

$R(LIN/\ff)$ can be generalized to a proof system
$R(PC_d/\ff)$, resolution over degree $d$ PC, operating with 
clauses formed by degree $\le d $ polynomials over $\ff$;
just add an extra rule
$$
\frac{C, g}{C, gh + h + 1}
$$
corresponding to the multiplication rule of polynomial calculus PC (cf.
Clegg, Edmonds and Impagliazzo \cite{CEI}).
Both processes from Section \ref{4.9.15c}
of reducing the width of clauses in a proof
work analogously as for $R(LIN/\ff)$. For definiteness let us now consider the construction
underlying Lemma \ref{4.4.15a}.
The clauses
$C^{\bf r} = \{g_1, \dots, g_w\}$ can be themselves replaced
by a single polynomial $1 - \Pi_{j \le w} (1-g_j)$
of degree $\le w d$. Hence the process can be repeated
any fixed number of times and thus, in fact, it can
be applied to $AC^0[2]$-formulas
and $AC^0[2]$-Frege proofs instead of $R(PC_d/\ff)$-proofs only.
This would result in a semantic PC-refutation of the original
set of clauses augmented by additional initial polynomials (analogous to axioms
$Ax(\pi,{\bf r})$) of degree $w^{O(1)}$ which yields also a syntactic 
PC-refutation of the same set of clauses and of the same degree
by Buss et.al.\cite[Thm.2.6]{BIKPRS}. A similar reduction can be obtained also
by using the characterization
of $AC^0[2]$-Frege proofs via the so called extended Nullstellensatz
proofs of Buss et.al.\cite{BIKPRS} and 
removing the extension axioms there by a random
assignment to the extension variables 
at the expense of introducing the new initial polynomials.
However, if monotone CLOs separating $\clique$ and $\color$ from Section \ref{20.3.16d} must be indeed large,
randomized feasible interpolation will not work in this situation as 
constant depth Frege systems admit short proofs of the weak pigeonhole principle and hence also
of the disjointness of the sets $\clique$ and $\color$ (when $\omega \geq 2 \xi$). 
Note also 
that $R(PC_d/\ff)$ even without the extra axioms p-simulates $R(d)$, a proof systems
operating with $d$-DNFs (cf. \cite{Kra-wphp}), which is known to be fairly strong (it corresponds
to bounded arithmetic theory $T^2_2(\alpha)$ for $d$ poly-logarithmic in $n$, cf.\cite{Kra-wphp}).

\bigskip
\noindent
{\bf Acknowledgements:} 

I thank Michal Garl\'{\i}k for pointing out a missing $\log n$ factor in Section \ref{4.9.15c},
to
Igor C. Oliveira and Pavel Pudl\' ak for comments on drafts of a part
of the paper and to Neil Thapen for discussions about related topics.

\bigskip
\noindent
{\bf Mailing address:}

Department of Algebra

Faculty of Mathematics and Physics

Charles University

Sokolovsk\' a 83, Prague 8, CZ - 186 75

The Czech Republic

{\tt krajicek@karlin.mff.cuni.cz}

\end{document}